\documentclass[english,11pt]{article}
\usepackage[latin1]{inputenc}
\usepackage{amsmath,amsthm,amssymb,babel}

\newtheorem{teo}{Theorem}

\newtheorem{prop}{Proposition}

\newtheorem{lemma}{Lemma}
\newtheorem{rem}{Remark}

\def\proof{{\it Proof.}\ }
\def\endproof{\hfill $\Box$\par\vskip3mm}

\def\eq#1{(\ref{#1})}

\def\neweq#1{\begin{equation}\label{#1}}
\def\endeq{\end{equation}}

\def\phi{\varphi}
\def\RR{{\mathbb R} }

\def\ZZ{{\mathbb Z} }

\def\di{\displaystyle}

\date{}

\title{\sc Existence and multiplicity of solutions for 
quasilinear nonhomogeneous problems: an Orlicz-Sobolev
space setting\thanks{
Correspondence address: Vicen\c{t}iu R\u{a}dulescu, Department of
Mathematics, University of Craiova,  200585 Craiova, Romania. E-mail: 
{\tt
vicentiu.radulescu@math.cnrs.fr}}}
\author{\sc Mihai Mih\u ailescu and Vicen\c{t}iu R\u{a}dulescu\\
\small
Department of Mathematics, University of Craiova,  200585 Craiova, 
Romania\\
\small
E-mail addresses: {\tt mmihailes@yahoo.com}\qquad {\tt 
vicentiu.radulescu@math.cnrs.fr}}

\begin{document}
\baselineskip16pt
\maketitle\noindent{\small{\sc Abstract}.
We study the boundary value problem $-{\rm div}(\log(1+
|\nabla u|^q)|\nabla u|^{p-2}\nabla u)=f(u)$ in $\Omega$,
$u=0$ on $\partial\Omega$, where $\Omega$ is a bounded domain in
$\RR^N$ with smooth boundary. We distinguish the cases 
where either $f(u)=-\lambda|u|^{p-2}u+|u|^{r-2}u$ or $f(u)=\lambda|u|^{p-2}u-|u|^{r-2}u$, with 
$p$, $q>1$ , $p+q<\min\{N,r\}$, and $r<(Np-N+p)/(N-p)$. 
In the first case we show the existence of infinitely many weak 
solutions for any $\lambda>0$. In the second case we prove the existence of a nontrivial weak solution if
$\lambda$ is sufficiently large. Our approach relies on adequate variational methods in Orlicz-Sobolev 
spaces.  \\
\small{\bf 2000 Mathematics
Subject Classification:}  35D05, 35J60, 35J70, 46N20. \\
\small{\bf Key words:}  nonhomogeneous operator,  Orlicz-Sobolev 
space, critical point, weak solution.}

\section{Introduction and main results}
Classical Sobolev and Orlicz-Sobolev spaces play a significant role in many 
fields of mathematics, such as approximation theory, partial differential equations, 
calculus of variations, non-linear potential theory, the theory of quasiconformal mappings, 
differential geometry, geometric function theory, and probability theory.
These spaces consists of functions that have weak derivatives and satisfy certain integrability conditions. 
The study of nonlinear elliptic equations involving quasilinear
homogeneous type operators is based on the theory of Sobolev spaces 
$W^{m,p}(\Omega)$ in order to find weak solutions. In the case of nonhomogeneous 
differential operators, the natural setting for this approach is
the use of Orlicz-Sobolev spaces. 
The basic idea is to replace the Lebesgue spaces 
$L^p(\Omega)$ by more general spaces $L_\Phi(\Omega)$, called {\it
Orlicz spaces}. The spaces $L_\Phi(\Omega)$ were thoroughly studied 
in the monograph by Kranosel'skii and Rutickii \cite{KR} and also
in the doctoral thesis of Luxemburg \cite{L}. If the role played
by $L^p(\Omega)$ in the definition of the Sobolev spaces $W^{m,p}
(\Omega)$ is assigned instead to an Orlicz space $L_\Phi(\Omega)$
the resulting space is denoted by $W^mL_\Phi(\Omega)$ and called
an {\it Orlicz-Sobolev space}. Many properties of Sobolev spaces
have been extended to Orlicz-Sobolev spaces, mainly by Dankert \cite{dank}, Donaldson
and Trudinger \cite{DT}, and O'Neill \cite{onei} (see also Adams \cite{A}
for an excellent account of those works). 

This paper is devoted to the study of weak solutions for problems 
of the type
\begin{equation}\label{Pr1}
\left\{\begin{array}{lll}
-{\rm div}(a(|\nabla u(x)|)\nabla u(x))=f(u(x)), 
&\mbox{for}& x\in\Omega\\
u(x)=0, &\mbox{for}& x\in\partial\Omega
\end{array}\right.
\end{equation}   
where $\Omega\subset\RR^N$ ($N\geq 3$) is a bounded domain with 
smooth boundary.

The first general existence result  using the theory of monotone operators in Orlicz-Sobolev spaces
were obtained in \cite{dona1} and in \cite{G, G1}. Other recent work that puts the problem into this framework
is contained in \cite{Clem1,Clem2,Gar,les}. In these papers, the existence results are obtained
using variational techniques, monotone operator methods or fixed point and degree theory arguments.

The case where $a(t)=t^{p-2}$ ($p>1$, $t\geq 0$) is fairly understood and a great variety of existence results are 
available. 
In this paper we focus on the case where
 $a:[0,\infty)
\rightarrow\RR$ is defined by $a(t)=\log(1+t^q)\cdot t^p$, where $p$, $q>1$. We treat
separately the cases where either
 $f(t)=-\lambda|t|^{p-2}t+|t|^{r-2}t$ or $f(t)=\lambda|t|^{p-2}t-|t|^{r-2}t$, where 
 $r<(Np-N+p)/(N-p)$ and $\lambda$ is a positive parameter.

We remark that we deal with a nonhomogeneous operator in the 
divergence form. Thus, we introduce an Orlicz-Sobolev 
space setting for problems of type \eq{Pr1}.

Define
$$\phi(t):=\log(1+|t|^q)\cdot|t|^{p-2}t,\qquad \mbox{for all}\ t\in\RR$$
and 
$$\Phi(t):=\int_0^t\phi(s),\qquad \mbox{for all}\ t\in\RR.$$
A straightforward computation yields 
$$\Phi(t)=\frac{1}{p}\log(1+|t|^q)\cdot|t|^p-\frac{q}{p}
\int_0^{|t|}\frac{s^{p+q-1}}{1+s^q}\;ds,$$
for all $t\in\RR.$
We point out that $\phi$ is an odd, increasing 
homeomorphism of $\RR$ into $\RR$, while $\Phi$ is convex and even
on $\RR$ and increasing from $\RR_+$ to $\RR_+$.

Set 
$$\Phi^\star(t):=\int_0^t\phi^{-1}(s)\;ds,\qquad\mbox{for all}\ t\in\RR.$$
The functions $\Phi$ and $\Phi^\star$ are complementary $N$-functions
(see \cite{A,KR,KJF}).

Define the Orlicz class
$$K_\Phi(\Omega):=\{u:\Omega\rightarrow\RR,\;{\rm measurable};\;
\int_\Omega\Phi(|u(x)|)\;dx<\infty\}$$
and the Orlicz space
$$L_\Phi(\Omega):=\;{\rm the}\;{\rm linear}\;{\rm hull}\;{\rm of}\;
K_\Phi(\Omega).$$
The space $L_\Phi(\Omega)$ is a Banach space endowed with the 
Luxemburg norm
$$\|u\|_\Phi:=\inf\left\{k>0;\ \int_\Omega\Phi\left(\frac{u(x)}{k}
\right)\;dx\leq 1\right\}$$
or the equivalent norm (the Orlicz norm)
$$\|u\|_{(\Phi)}:=\sup\left\{\left|\int_\Omega uvdx\right|;\ v\in K_{\overline\Phi}(\Omega),\ 
\int_\Omega \overline\Phi (|v|)dx\leq 1\right\}\,,$$
where $\overline\Phi$ denotes the conjugate Young function of $\Phi$, that is,
$$\overline \Phi (t)=\sup\{ts-\Phi (s);\ s\in\RR\}\,.$$

By Lemma 2.4 and Example 2 in \cite[p.~243]{Clem2} we have
$$1<\liminf\limits_{t\rightarrow\infty}\frac{t\phi(t)}{\Phi(t)}
\leq\sup\limits_{t>0}\frac{t\phi(t)}{\Phi(t)}<\infty.$$  
The above inequalities imply that $\Phi$ satisfies the 
$\Delta_2$-condition. By Lemma C.4 in \cite{Clem2} it follows that
$\Phi^\star$ also satisfies the $\Delta_2$-condition. 
Then, according to \cite{A}, p. 234, it
folows that $L_\Phi(\Omega)=K_\Phi(\Omega)$. Moreover, by 
Theorem 8.19 in \cite{A} $L_\Phi(\Omega)$ is reflexive.

We denote by $W^1L_\Phi(\Omega)$ the Orlicz-Sobolev space defined
by
$$W^1L_\Phi(\Omega):=\left\{u\in L_\Phi(\Omega);\;\frac{\partial u}
{\partial x_i}\in L_\Phi(\Omega),\;i=1,...,N\right\}.$$
This is a Banach space with respect to the norm
$$\|u\|_{1,\Phi}:=\|u\|_\Phi+\||\nabla u|\|_\Phi.$$
We also define the Orlicz-Sobolev space $W_0^1L_\Phi(\Omega)$ as
the closure of $C_0^\infty(\Omega)$ in $W^1L_\Phi(\Omega)$. By
Lemma 5.7 in \cite{G} we obtain that on $W_0^1L_\Phi(\Omega)$ we
may consider an equivalent norm
$$\|u\|:=\||\nabla u|\|_\Phi.$$ 
The space $W_0^1L_\Phi(\Omega)$ is also a reflexive Banach space.

\medskip
In the first part of the present paper we study the boundary value problem
\begin{equation}\label{1}
\left\{\begin{array}{lll}
-{\rm div}(\log(1+|\nabla u(x)|^q)|\nabla u(x)|^{p-2}\nabla 
u(x))=-\lambda|u(x)|^{p-2}u(x)+|u(x)|^{r-2}u(x), 
&\mbox{for}& x\in\Omega\\
u(x)=0, &\mbox{for}& x\in\partial\Omega.
\end{array}\right.
\end{equation}  
We say that $u\in W_0^1L_\Phi(\Omega)$ is a {\it weak solution} of
problem \eq{1} if 
$$\int_\Omega\log(1+|\nabla u(x)|^q)|\nabla u(x)|^{p-2}
\nabla u\nabla v\;dx+\lambda\int_\Omega|u(x)|^{p-2}u(x)v(x)\;dx
-\int_\Omega|u(x)|^{r-2}u(x)v(x)\;dx=0$$
for all $v\in W_0^1L_\Phi(\Omega)$.

We prove the following multiplicity result.
\begin{teo}\label{t1}
Assume that $p$, $q>1$ , $p+q<N$, $p+q<r$ and $r<(Np-N+p)/(N-p)$. Then
for every $\lambda>0$ problem \eq{1} has infinitely many weak
solutions.
\end{teo}

Next, we consider the problem
\begin{equation}\label{2}
\left\{\begin{array}{lll}
-{\rm div}(\log(1+|\nabla u(x)|^q)|\nabla u(x)|^{p-2}\nabla 
u(x))=\lambda|u(x)|^{p-2}u(x)-|u(x)|^{r-2}u(x), 
&\mbox{for}& x\in\Omega\\
u(x)=0, &\mbox{for}& x\in\partial\Omega.
\end{array}\right.
\end{equation}  
We say that $u\in W_0^1L_\Phi(\Omega)$ is a {\it weak solution} of
problem \eq{2} if 
$$\int_\Omega\log(1+|\nabla u(x)|^q)|\nabla u(x)|^{p-2}
\nabla u\nabla v\;dx-\lambda\int_\Omega|u(x)|^{p-2}u(x)v(x)\;dx
+\int_\Omega|u(x)|^{r-2}u(x)v(x)\;dx=0$$
for all $v\in W_0^1L_\Phi(\Omega)$.

We prove
\begin{teo}\label{t2}
Assume that the hypotheses of Theorem \ref{t1} are fulfilled. Then
there exists $\lambda_\star>0$ such that for any $\lambda\geq
\lambda_\star$, problem \eq{2} has a nontrivial weak solution.
\end{teo}

\section{Auxiliary results on Orlicz-Sobolev embeddings}
In many applications of Orlicz-Sobolev spaces to boundary value problems for 
nonlinear partial differential equations, the compactness of the embeddings plays a central role.
Compact embedding theorems for Sobolev or Orlicz-Sobolev spaces are also intimately connected with the problem
of discreteness of spectra of Schr\"odinger operators (see Benci and Fortunato \cite{benci} and Reed and Simon \cite{reed}).

While the Banach spaces $W^1L_\Phi(\Omega)$ and $W_0^1L_\Phi(\Omega)$ 
can be defined from fairly general convex properties of $\Phi$,
it is also well known that the specific functional-analytic and topological properties of
these spaces depend very sensitively on the rate of growth of $\Phi$ at infinity. Compactness is
not an exception and, using standard notions traditionally used to describe convex functions, 
we recall in this section a compact embedding theorem for a class of
Orlicz-Sobolev spaces. 

Define 
the Orlicz-Sobolev conjugate $\Phi_\star$ of $\Phi$ by
$$\Phi_\star^{-1}(t):=\int_0^t\frac{\Phi^{-1}(s)}
{s^{\frac{N+1}{N}}}\;ds.$$

\begin{prop}\label{p1}
Assume that the hypotheses of Theorems \ref{t1} or \ref{t2} are fulfilled. Then the following
properties hold true.
\smallskip

\noindent a) $\di\lim_{t\rightarrow 0}\di\int_t^1\di
\frac{\Phi^{-1}(s)}{s^{\frac{N+1}{N}}}\;ds<\infty$;
\smallskip

\noindent b) $\di\lim_{t\rightarrow\infty}\di\int_1^t
\di\frac{\Phi^{-1}(s)}{s^{\frac{N+1}{N}}}\;ds=\infty$;
\smallskip

\noindent c) $\di\lim_{t\rightarrow\infty}\di\frac{|t|^{\gamma+1}}
{\Phi_\star(kt)}=0$, for all $k>0$ and all $1\leq\gamma<
\di\frac{Np-N+p}{N-p}$.  
\end{prop}
\proof
a) By L'H\^opital's rule we have
\begin{eqnarray*}
\lim_{t\searrow 0}\frac{\Phi(t)}{t^{p+q}}&=&\lim_{t\searrow 0}
\frac{\phi(t)}{(p+q)t^{p+q-1}}\\
&=&\frac{1}{p+q}\lim_{t\searrow 0}\frac{\log(1+t^q)}{t^q}=
\frac{1}{p+q}\lim_{t\searrow 0}\frac{\frac{qt^{q-1}}{1+t^q}}
{qt^{q-1}}=\frac{1}{p+q}.
\end{eqnarray*}
We deduce that $\Phi$ is equivalent to $t^{p+q}$ near zero. Using
that fact and the remarks on p. 248 in \cite{A} we infer that a)
holds true if and only if
$$\lim_{t\rightarrow 0}\int_t^1
\frac{s^{\frac{1}{p+q}}}{s^{\frac{N+1}{N}}}\;ds<\infty,$$
or
$$p+q<N.$$
The last inequality holds since the hypotheses of Theorems \ref{t1}
or \ref{t2} are fulfilled. 

b) By the change of variable $s=\Phi(\tau)$ we obtain
\begin{equation}\label{e0}
\int_1^t\frac{\Phi^{-1}(s)}{s^{\frac{N+1}{N}}}\;ds=
\int_{\Phi^{-1}(1)}^{\Phi^{-1}(t)}
\frac{\tau\phi(\tau)}{\Phi(\tau)}(\Phi(\tau))^{-1/N}d\tau.
\end{equation}
A simple calculation yields
$$0\leq\lim_{\tau\rightarrow\infty}
\frac{\di\int_0^\tau\frac{s^{p+q-1}}{1+s^q}\;ds}
{\tau^p\log(1+\tau^q)}\leq\lim_{\tau\rightarrow\infty}
\frac{\di\int_0^\tau\frac{s^{p+q-1}}{s^q}\;ds}{\tau^p\log(1+\tau^q)}
=\lim_{\tau\rightarrow\infty}\frac{\di\frac{1}{p}\tau^p}{\tau^p
\log(1+\tau^q)}=0.$$
Thus
\begin{equation}\label{e1}
\lim_{\tau\rightarrow\infty}\frac{\di\int_0^\tau
\frac{s^{p+q-1}}{1+s^q}\;ds}{\tau^p\log(1+\tau^q)}=0.
\end{equation}
A first consequence of the above relation is that
\begin{equation}\label{e2}
\lim_{t\rightarrow\infty}\frac{\Phi(t)}{t^p\log(1+t^q)}=\frac{1}{p}.
\end{equation}
On the other hand, by \eq{e1},
\begin{equation}\label{e3}
\begin{array}{lll}
\lim\limits_{\tau\rightarrow\infty}\di\frac{\tau\phi(\tau)}
{\Phi(\tau)}
&=&\lim\limits_{\tau\rightarrow\infty}\di\frac{\tau^p\log(1+\tau^q)}
{\di\frac{1}{p}\tau^p\log(1+\tau^q)-\di\frac{q}{p}\di\int_0^\tau\frac
{s^{p+q-1}}{1+s^q}\;ds}\\
&=&p\lim\limits_{\tau\rightarrow\infty}
\left(1-q\cdot\di\frac{\di\int_0^\tau\di\frac
{s^{p+q-1}}{1+s^q}\;ds}{\tau^p\log(1+\tau^q)}\right)^{-1}=p
\end{array}
\end{equation}
and
\begin{equation}\label{e4}
\lim_{t\rightarrow\infty}\Phi(t)=\lim_{t\rightarrow\infty}
\frac{1}{p}t^p\log(1+t^q)\left[1-q\cdot\di\frac{\di\int_0^{|t|}
\di\frac{s^{p+q-1}}{1+s^q}\;ds}{t^p\log(1+t^q)}\right]=\infty.
\end{equation}
Relations \eq{e0}, \eq{e3} and \eq{e4} yield
$$\lim_{t\rightarrow\infty}\int_1^t
\frac{\Phi^{-1}(s)}{s^{\frac{N+1}{N}}}\;ds=\infty\,.$$
Equivalently, we can write
$$\int_{\Phi^{-1}(1)}^\infty\frac{d\tau}{[\Phi(\tau)]^{1/N}}=
\infty$$
or, by \eq{e2},
\begin{equation}\label{e5}
\int_{\Phi^{-1}(1)}^\infty\frac{d\tau}{\tau^{p/N}[\log(1+\tau^q)]
^{1/N}}=\infty.
\end{equation}
Since
$$\log(1+\theta)\leq\theta,\;\;\;\forall\;\theta>0$$
we deduce that
$$\frac{1}{{\tau^{p/N}[\log(1+\tau^q)]^{1/N}}}\geq
\frac{1}{\tau^{(p+q)/N}},\;\;\;\forall\;\tau>0.$$
Since
$p+q<N$, we find
$$\int_{\Phi^{-1}(1)}^\infty\tau^{-(p+q)/N}d\tau=\infty$$
and thus relation \eq{e5} holds true. We conclude that
$$\lim_{t\rightarrow\infty}\int_1^t
\frac{\Phi^{-1}(s)}{s^{\frac{N+1}{N}}}\;ds=\infty$$

c) Let $\gamma$ be fixed such that $1\leq\gamma<(Np-N+p)/(N-p)$.

By Adams \cite[p. 231]{A}, we have 
$$\lim_{t\rightarrow\infty}\frac{|t|^{\gamma+1}}
{\Phi_\star(kt)}=0,\;\;\; \forall\;k>0$$
if and only if
\begin{equation}\label{e6}
\lim_{t\rightarrow\infty}
\frac{\Phi_\star^{-1}(t)}{t^{1/(\gamma+1)}}=0.
\end{equation}
Using again L'H\^opital's rule we deduce that
$$\limsup_{t\rightarrow\infty}
\frac{\Phi_\star^{-1}(t)}{t^{1/(\gamma+1)}}\leq
(\gamma+1)\limsup_{t\rightarrow\infty}\frac{\Phi^{-1}(t)}
{t^{\frac{1}{\gamma+1}+\frac{1}{N}}}.$$
Setting $\tau=\Phi(t)$ we obtain
$$\limsup_{t\rightarrow\infty}
\frac{\Phi_\star^{-1}(t)}{t^{1/(\gamma+1)}}\leq(\gamma+1)
\limsup_{\tau\rightarrow\infty}\frac{\tau}{[\Phi(\tau)]^
{\frac{1}{\gamma+1}+\frac{1}{N}}}.$$
Since $\gamma<(Np-N+p)/(N-p)$ we have
$$p>\frac{N(\gamma+1)}{N+\gamma+1}.$$
Using the above inequality and \eq{e1} we get
$$\limsup_{\tau\rightarrow\infty}\frac{\tau^
{\frac{N(\gamma+1)}{N+\gamma+1}}}{\Phi(\tau)}=0.$$
We conclude that c) holds true.

Thus the proof of Proposition \ref{p1} is complete.      \endproof

\begin{rem}\label{r1}
Proposition \ref{p1} enables us to apply Theorem 2.2 in \cite{Gar}
(see also Theorem 8.33 in \cite{A}) in order to obtain that
$W_0^1L_\Phi(\Omega)$ is compactly embedded in $L^{\gamma+1}
(\Omega)$ provided that $1\leq\gamma<(Np-N+p)/(N-p)$.
\end{rem}

An important role in what follows will be played by
$$p^0:=\sup_{t>0}\frac{t\phi(t)}{\Phi(t)}.$$
\begin{rem}\label{r2}
By Example 2 on p. 243 in \cite{Clem2} it follows that
$$p^0=p+q.$$
\end{rem}

\section{Proof of Theorem \ref{t1}}
The key argument in the proof of Theorem \ref{t1} is the following
$\ZZ_2$-symmetric version (for even functionals) of the Mountain 
Pass Lemma (see Theorem 9.12 in \cite{R}).

\smallskip
\noindent{\bf Mountain Pass Lemma.} {\it Let $X$ be an infinite 
dimensional real Banach space and let $I\in C^1(X,\RR)$ be even, 
satisfying the Palais-Smale condition (that is, any sequence $\{x_n\}
\subset X$ such that $\{I(x_n)\}$ is bounded and $I^{'}(x_n)
\rightarrow 0$ in $X^\star$ has a convergent subsequence) and 
$I(0)=0$. Suppose that
\smallskip

\noindent (I1) There exist two constants $\rho$, $b>0$ such that
$I(x)\geq b$ if $\|x\|=\rho.$
\smallskip

\noindent (I2) For each finite dimensional subspace $X_1\subset X$,
the set $\{x\in X_1;\;I(x)\geq 0\}$ is bounded.
\smallskip

Then $I$ has an unbounded sequence of critical values.} 

Let $E$ denote the Orlicz-Sobolev space $W_0^1L_\Phi(\Omega)$.
Let $\lambda>0$ be arbitrary but fixed.

The energy functional associated to problem \eq{1} is
$J_\lambda:E\rightarrow\RR$ defined by
$$J_\lambda(u):=\int_\Omega\Phi(|\nabla u(x)|)\;dx+\frac{\lambda}{p}
\int_\Omega|u(x)|^p\;dx-\frac{1}{r}\int_\Omega|u(x)|^r\;dx.$$
By Remark \ref{r1}, $J_\lambda$ is well defined
on $E$.

Let us denote by $J_{\lambda,1}$, $J_{\lambda,2}:E\rightarrow\RR$
the functionals
$$J_{\lambda,1}(u):=\int_\Omega\Phi(|\nabla u(x)|)\;dx\;\;\;
{\rm and}\;\;\;J_{\lambda,2}(u):=\frac{\lambda}{p}
\int_\Omega|u(x)|^p\;dx-\frac{1}{r}\int_\Omega|u(x)|^r\;dx.$$
Therefore
$$J_\lambda(u)=J_{\lambda,1}(u)+J_{\lambda,2}(u),\;\;\;\forall\;
u\in E.$$
By Lemma 3.4 in \cite{Gar} it follows that $J_{\lambda,1}$ is a
$C^1$ functional, with the Fr\'echet derivative given by
$$\langle J_{\lambda,1}^{'}(u),v\rangle=\int_\Omega\log(1+|\nabla 
u(x)|^q)|\nabla u(x)|^{p-2}\nabla u(x)\nabla v(x)\;dx\,,$$
for all $u$, $v\in E$.

Similar arguments as those used in the proof of Lemma 2.1 in
\cite{Clem1} imply that $J_{\lambda,2}$ is of class $C^1$ with
the Fr\'echet derivative given by
$$\langle J_{\lambda,2}^{'}(u),v\rangle=\lambda\int_\Omega
|u(x)|^{p-2}u(x)v(x)\;dx-\int_\Omega|u(x)|^{r-2}u(x)v(x)\;dx\,,$$
for all $u$, $v\in E$.

The above information shows that $J_\lambda\in C^1(E,\RR)$ and
\begin{eqnarray*}
\langle J_\lambda^{'}(u),v\rangle&=&\int_\Omega\log(1+|\nabla 
u(x)|^q)|\nabla u(x)|^{p-2}\nabla u(x)\nabla v(x)\;dx\\
&+&\lambda\int_\Omega|u(x)|^{p-2}u(x)v(x)\;dx-\int_\Omega
|u(x)|^{r-2}u(x)v(x)\;dx
\end{eqnarray*}    
for all $u$, $v\in E$. Thus, the weak solutions of
\eq{1} coincide with the critical points of $J_\lambda$. 

\begin{lemma}\label{le1}
There exist $\eta>0$ and $\alpha>0$ such that 
$J_\lambda(u)\geq\alpha>0$ for any $u\in E$ with $\|u\|=\eta$.
\end{lemma}
\proof
In order to prove Lemma \ref{le1} we first show that
\begin{equation}\label{l2}
\Phi(t)\geq\tau^{p^0}\Phi(t/\tau),\;\;\;\forall\; t>0\;{\rm and}\;
\tau\in(0,1]\,,
\end{equation}
where $p^0$ is defined in the previous section.

Indeed, since
$$p^0=\sup_{t>0}\frac{t\phi(t)}{\Phi(t)}$$
we have
$$\frac{t\phi(t)}{\Phi(t)}\leq p^0,\;\;\;\forall\;t>0.$$
Let $\tau\in(0,1]$ be fixed. We have
$$\log(\Phi(t/\tau))-\log(\Phi(t))=\int_t^{t/\tau}\frac{\phi(s)}
{\Phi(s)}\;ds\leq\int_t^{t/\tau}\frac{p^0}{s}\;ds=\log(\tau^{-p^0})$$
and it follows that \eq{l2} holds true.

Fix $u\in E$ with $\|u\|<1$ and $\xi\in(0,\|u\|)$. Using relation
\eq{l2} we have
\begin{equation}\label{for1}
\int_\Omega\Phi(|\nabla u(x)|)\;dx\geq\xi^{p^0}\int_\Omega
\Phi\left(\frac{|\nabla u(x)|}{\xi}\right)\;dx.
\end{equation}
Defining $v(x)=|\nabla u(x)|/\xi$, for all $x\in\Omega$, we have
$\|v\|_{\Phi}=\|u\|/\xi>1$. Since $\Phi(t)\leq\frac{t\phi(t)}
{p}$, for all $t\in\RR$, by Lemma C.9 in \cite{Clem2} we deduce that
\begin{equation}\label{for2}
\int_\Omega\Phi(v(x))\;dx\geq\|v\|_{\Phi}^p>1.
\end{equation}
Relations \eq{for1} and \eq{for2} show that
$$\int_\Omega\Phi(|\nabla u(x)|)\;dx\geq\xi^{p^0}.$$
Letting $\xi\nearrow\|u\|$ in the above inequality we obtain
\begin{equation}\label{l3}
\int_\Omega\Phi(|\nabla u(x)|)\;dx\geq\|u\|^{p^0},\;\;\;\forall\;
u\in E\;{\rm with}\;\|u\|<1.
\end{equation}
On the other hand, since $E$ is continuously embedded in 
$L^r(\Omega)$, it follows that there exists a positive constant
$C_1>0$ such that
\begin{equation}\label{l4}
\int_\Omega|u(x)|^r\;dx\leq C_1\cdot\|u\|^r,\;\;\;\forall\;u\in E.
\end{equation}
Using relations \eq{l3} and \eq{l4} we deduce that for all $u\in E$
with $\|u\|\leq 1$ we have
\begin{eqnarray*}
J_\lambda(u)&\geq&\int_\Omega\Phi(|\nabla u(u)|)\;dx-\frac{1}{r}
\int_\Omega|u(x)|^r\;dx\\
&\geq&\|u\|^{p^0}-\frac{C_1}{r}\cdot\|u\|^r\\
&=&\left(1-\frac{C_1}{r}\cdot\|u\|^{r-p^0}\right)\|u\|^{p^0}.
\end{eqnarray*}
But, by Remark \ref{r2} and the hypotheses of Theorem \ref{t1}, 
we have  $p^0=p+q<r$. We conclude that Lemma \ref{le1} holds
true.    \endproof

\begin{lemma}\label{le2}
Assume that $E_1$ is a finite dimensional subspace of $E$. Then the set $S=\{
u\in E_1;\;J_\lambda(u)\geq 0\}$ is bounded.
\end{lemma}
\proof
With the same arguments as those used in the proof of relation 
\eq{l2} we have 
\begin{equation}\label{l41}
\frac{\Phi(\sigma t)}{\Phi(t)}\leq\sigma^{p^0},\;\;\;\forall\;
t>0\;{\rm and}\;\sigma>1.
\end{equation}
Then, for all $u\in E$ with $\|u\|>1$, relation \eq{l41} implies
\begin{equation}\label{l5}
\begin{array}{lll}
\di\int_\Omega\Phi(|\nabla u(x)|)\;dx&=&\di\int_\Omega\Phi\left(
\|u\|\di\frac{|\nabla u(x)|}{\|u\|}\right)\;dx\\
&\leq&\|u\|^{p^0}\di\int_\Omega\Phi\left(\di\frac{|\nabla u(x)|}
{\|u\|}\right)\;dx\\
&\leq&\|u\|^{p^0}.
\end{array}
\end{equation}
On the other hand, since $E$ is continuously embedded in 
$L^p(\Omega)$ it follows that there exists a positive constant
$C_2>0$ such that
\begin{equation}\label{l6}
\int_\Omega|u(x)|^p\;dx\leq C_2\cdot\|u\|^p,\;\;\;\forall\;u\in E.
\end{equation}
Relations \eq{l5} and \eq{l6} yield
\begin{equation}\label{l7}
J_\lambda(u)\leq\|u\|^{p^0}+\frac{\lambda}{p}\cdot C_2\cdot\|u\|^p
-\frac{1}{r}\int_\Omega|u(x)|^r\;dx\,,
\end{equation}
for all $u\in E$ with $\|u\|>1$.

We point out that the functional $|\cdot|_{r}:E\rightarrow\RR$ 
defined by
$$|u|_{r}=\left(\int_\Omega|u(x)|^{r}\;dx\right)^{1/{r}}$$
is a norm in $E$. In the finite dimensional subspace $E_1$ the norms
$|.|_{r}$ and $\|.\|$ are equivalent, so there exists a
positive constant $C_3=C_3(E_1)$ such that
$$\|u\|\leq C_3\cdot|u|_{r},\;\;\;\forall\;u\in E_1.$$
The above remark and relation \eq{l7} imply
$$J_\lambda(u)\leq\|u\|^{p^0}+\frac{\lambda}{p}\cdot C_2\cdot
\|u\|^p-\frac{1}{r}\cdot C_3^{-1}\cdot\|u\|^r\,,$$
for all $u\in E_1$ with $\|u\|>1$.

Hence
\begin{equation}\label{l8}
\|u\|^{p^0}+\frac{\lambda}{p}\cdot C_2\cdot\|u\|^p-\frac{1}{r}
\cdot C_3^{-1}\cdot\|u\|^r\geq 0\,,
\end{equation}
for all $u\in S$ with $\|u\|>1$. Since, by Remark \ref{r2} and
the hypotheses of Theorem \ref{t1} we have $r>p^0>p$, the above
relation implies that $S$ is bounded in $E$.  \endproof

\begin{lemma}\label{le3}
Assume that $\{u_n\}\subset E$ is a sequence which satisfies the properties
\begin{equation}\label{l9}
|J_\lambda(u_n)|<M
\end{equation}
\begin{equation}\label{l10}
J_\lambda^{'}(u_n)\rightarrow 0\;\;\;{\rm as}\;n\rightarrow\infty\,,
\end{equation}
where $M$ is a positive constant. Then $\{u_n\}$ possesses a 
convergent subsequence.
\end{lemma}
\proof
First, we show that $\{u_n\}$ is bounded in $E$. Assume by 
contradiction the contrary. Then, passing eventually to a subsequence,
still denoted by $\{u_n\}$, we may assume that $\|u_n\|
\rightarrow\infty$ as $n\rightarrow\infty$. Thus we may consider that
$\|u_n\|>1$ for any integer $n$.

By \eq{l10} we deduce that there exists $N_1>0$ such that for any
$n>N_1$ we have
$$\|J_\lambda^{'}(u_n)\|\leq 1.$$
On the other hand, for any $n>N_1$ fixed, the application
$$E\ni v\rightarrow\langle J_\lambda^{'}(u_n),v\rangle$$
is linear and continuous.

The above information yields 
$$|\langle J_\lambda^{'}(u_n),v\rangle|\leq
\|J_\lambda^{'}(u_n)\|\cdot\|v\|\leq\|v\|,\;\;\;\forall
v\in E,\;n>N_1.$$
Setting $v=u_n$ we have
$$-\|u_n\|\leq\int_\Omega\log(1+|\nabla u_n(u)|^q)|\nabla u_n(x)|^p
\;dx+\lambda\int_\Omega|u_n(x)|^{p}\;dx-
\int_\Omega|u_n(x)|^{r}\;dx\leq\|u_n\|,$$
for all $n>N_1$. We obtain
\begin{equation}\label{l11}
-\|u_n\|-\int_\Omega\log(1+|\nabla u_n(u)|^q)|\nabla u_n(x)|^p
\;dx-\lambda\int_\Omega|u_n(x)|^{p}\;dx\leq-\int_\Omega|u_n(x)|^{r}
\;dx,
\end{equation}
for any $n>N_1$. 

If $\|u_n\|>1$, then relations \eq{l9} and \eq{l11}  imply
\begin{eqnarray*}
M>J_\lambda(u_n)&=&\int_\Omega\Phi(|\nabla u_n(x)|)\;dx+
\frac{\lambda}{p}\int_\Omega|u_n(x)|^p\;dx-\frac{1}{r}
\int_\Omega|u_n(x)|^r\;dx\\
&\geq&\int_\Omega\Phi(|\nabla u_n(x)|)\;dx+\lambda\cdot\left(
\frac{1}{p}-\frac{1}{r}\right)\cdot\int_\Omega|u_n(x)|^p\;dx-\\
& &\frac{1}{r}\cdot\int_\Omega\log(1+|\nabla u_n(u)|^q)
|\nabla u_n(x)|^p\;dx-\frac{1}{r}\cdot\|u_n\|\\
&=&\int_\Omega\Phi(|\nabla u_n(x)|)\;dx-\frac{1}{r}\cdot
\int_\Omega\phi(|\nabla u_n(x)|)|\nabla u_n(x)|\;dx+\\
& &\lambda\cdot\left(\frac{1}{p}-\frac{1}{r}\right)\cdot
\int_\Omega|u_n(x)|^p\;dx-\frac{1}{r}\cdot\|u_n\|.
\end{eqnarray*}
Since 
$$p^0\geq\frac{t\phi(t)}{\Phi(t)},\;\;\;\forall\;t>0$$
we find 
$$\int_\Omega\Phi(|\nabla u_n(x)|)\;dx-\frac{1}{r}\cdot
\int_\Omega\phi(|\nabla u_n(x)|)|\nabla u_n(x)|\;dx\geq
\left(1-\frac{p^0}{r}\right)\int_\Omega\Phi(|\nabla u_n(x)|)\;dx.$$
Using the above relations we deduce that for any $n>N_1$ such that 
$\|u_n\|>1$ we have
\begin{equation}\label{l12}
M>\left(1-\frac{p^0}{r}\right)\cdot\int_\Omega\Phi(|\nabla u_n(x)|)
\;dx-\frac{1}{r}\cdot\|u_n\|.
\end{equation}

Since $\Phi(t)\leq(t\phi(t))/p$ for all $t\in\RR$ we deduce by
Lemma C.9 in \cite{Clem2} that
\begin{equation}\label{l13}
\int_\Omega\Phi(|\nabla u_n(x)|)\;dx\geq\|u_n\|^p,
\end{equation}
for all $n>N_1$ with $\|u_n\|>1$.

Relations \eq{l12} and \eq{l13} imply
$$M>\left(1-\frac{p^0}{r}\right)\cdot\|u_n\|^p-\frac{1}{r}\cdot
\|u_n\|,$$
for all $n>N_1$ with $\|u_n\|>1$. Since $p^0<r$, letting 
$n\rightarrow\infty$ we obtain a contradiction. It follows that
$\{u_n\}$ is bounded in $E$.

Since $\{u_n\}$ is bounded in $E$ we deduce that there exists a
subsequence, still denoted by $\{u_n\}$, and $u_0\in E$ such that
$\{u_n\}$ converges weakly to $u_0$ in $E$. Since $E$ is compactly
embedded in $L^p(\Omega)$ and $L^r(\Omega)$ it follows
that $\{u_n\}$ converges strongly to $u_0$ in $L^p(\Omega)$ 
and $L^r(\Omega)$. Hence
\begin{equation}\label{l14}
\lim_{n\rightarrow\infty}J_{\lambda,2}(u_n)=J_{\lambda,2}(u_0)
\;\;\;{\rm and}\;\;\;\lim_{n\rightarrow\infty}J_{\lambda,2}^{'}
(u_n)=J_{\lambda,2}^{'}(u_0).
\end{equation} 
Since
$$J_{\lambda,1}(u)=J_\lambda(u)-J_{\lambda,2}(u),\;\;\;\forall\;
u\in E\,,$$
relations \eq{l14} and \eq{l10} imply
\begin{equation}\label{l15}
\lim_{n\rightarrow\infty}J_{\lambda,1}^{'}(u_n)=
-J_{\lambda,2}^{'}(u_0),\;\;\;{\rm in}\;E^\star.
\end{equation}
Using the fact that $\Phi$ is convex and thus $J_{\lambda,1}$ is
convex we have that
$$J_{\lambda,1}(u_n)\leq J_{\lambda,1}(u_0)+
\langle J_{\lambda,1}^{'}(u_n),u_n-u_0\rangle.$$
Passing to the limit as $n\rightarrow\infty$ and using \eq{l15}
we deduce that
\begin{equation}\label{l16}
\limsup_{n\rightarrow\infty}J_{\lambda,1}(u_n)\leq 
J_{\lambda,1}(u_0).
\end{equation}
Using again the fact that $J_{\lambda,1}$ is convex, it follows that
$J_{\lambda,1}$ is weakly lower semicontinuous and hence
\begin{equation}\label{l17}
\liminf_{n\rightarrow\infty}J_{\lambda,1}(u_n)\geq 
J_{\lambda,1}(u_0).
\end{equation}
By \eq{l16} and \eq{l17} we find
$$\lim_{n\rightarrow\infty}J_{\lambda,1}(u_n)= 
J_{\lambda,1}(u_0)$$
or
\begin{equation}\label{l18}
\lim_{n\rightarrow\infty}\int_\Omega\Phi(|\nabla u_n(x)|)\;dx=
\int_\Omega\Phi(|\nabla u_0(x)|)\;dx\,.
\end{equation}
Since $\Phi$ is increasing and convex, it follows that
$$\Phi\left(\frac{1}{2}|\nabla u_n(x)-\nabla u_0(x)|\right)\leq
\Phi\left(\frac{1}{2}\left(|\nabla u_n(x)|+|\nabla u_0(x)|\right)
\right)\leq\frac{\Phi(|\nabla u_n(x)|)+\Phi(|\nabla u_0(x)|)}{2}\,,$$
for all $x\in\Omega$ and all $n$. Integrating the above inequalities 
over $\Omega$ we find
$$0\leq\int_\Omega\Phi\left(\frac{1}{2}|\nabla (u_n-u_0)(x)|
\right)\;dx\leq\frac{\di\int_\Omega\Phi(|\nabla u_n(x)|)\;dx+
\di\int_\Omega\Phi(|\nabla u_0(x)|)\;dx}{2}\,,$$
for all $n$. We point out that Lemma C.9 in \cite{Clem2} implies
$$\int_\Omega\Phi(|\nabla u_n(x)|)\;dx\leq\|u_n\|^p<1,\;\;\;
{\rm provided}\;{\rm that}\;\|u_n\|<1\,,$$
while relation \eq{l41} yields
$$\int_\Omega\Phi(|\nabla u_n(x)|)\;dx\leq\|u_n\|^{p^0},\;\;\;
{\rm provided}\;{\rm that}\;\|u_n\|>1.$$
Since $\{u_n\}$ is bounded in $E$, the above inequalities prove the
existence of a positive constant $K_1$ such that
$$\int_\Omega\Phi(|\nabla u_n(x)|)\;dx\leq K_1,$$
for all $n$.
So, there exists a positive constant $K_2$ such that
\begin{equation}\label{l19}
0\leq\int_\Omega\Phi\left(\frac{1}{2}|\nabla (u_n-u_0)(x)|
\right)\;dx\leq K_2,
\end{equation} 
for all $n$.

On the other hand, since $\{u_n\}$ converges weakly to $u_0$ in $E$,
Theorem 2.1 in \cite{Gar} implies 
$$\int_\Omega\frac{\partial u_n}{\partial x_i}v\;dx
\rightarrow\int_\Omega\frac{\partial u_0}{\partial x_i}v\;dx,
\;\;\;\forall\;v\in L_{\Phi^\star}(\Omega),\;
\forall\;i=1,...,N.$$
In particular this holds for all $v\in L^\infty(\Omega)$. Hence
$\{\frac{\partial u_n}{\partial x_i}\}$ converges weakly to
$\frac{\partial u_0}{\partial x_i}$ in $L^1(\Omega)$ for all
$i=1,...,N$. Thus we deduce that
\begin{equation}\label{l20}
\nabla u_n(x)\rightarrow\nabla u_0(x)\;\;\;{\rm a.e.}\;x\in\Omega.
\end{equation}
Relations \eq{l19} and \eq{l20} and Lebesgue's dominated convergence theorem imply
$$\lim_{n\rightarrow\infty}\int_\Omega
\Phi\left(\frac{1}{2}|\nabla (u_n-u_0)(x)|\right)\;dx=0.$$
Taking into account that $\Phi$ satisfies the $\Delta_2$-condition
it follows by Lemma A.4 in \cite{Clem2} (see also \cite{A}, p. 236)
that
$$\lim_{n\rightarrow\infty}\left\|\frac{1}{2}(u_n-u_0)\right\|=0$$
and thus
$$\lim_{n\rightarrow\infty}\|(u_n-u_0)\|=0.$$
The proof of Lemma \ref{le3} is complete.  \endproof

{\sc Proof of Theorem \ref{t1} completed.}
It is clear that the functional $J_\lambda$ is even and verifies
$J_\lambda(0)=0$. Lemma \ref{le3} implies that $J_\lambda$ 
satisfies the Palais-Smale condition. On the other hand, 
Lemmas \ref{le1} and \ref{le2} show that conditions (I1) and (I2)
are satisfied. Thus the Mountain Pass Lemma can be applied to the
functional $J_\lambda$. We conclude that equation \eq{1} has 
infinitely many weak solutions in $E$.
The proof of Theorem \ref{t1} is complete.  \endproof
\begin{rem}\label{r3}
We point out the fact that the Orlicz-Sobolev space $E$ cannot
be replaced by a classical Sobolev space, since, in this case, 
condition (I1) in the Mountain Pass Lemma cannot be satisfied.
For a proof of that fact one can consult the proof of Remark 4
in \cite{Clem1} (p. 56-57).
\end{rem}

\section{Proof of Theorem \ref{t2}}

Let $\lambda>0$ be arbitrary but fixed.
Let $I_\lambda:E\rightarrow\RR$ be defined by
$$I_\lambda(u):=\int_\Omega\Phi(|\nabla u(x)|)\;dx-\frac{\lambda}{p}
\int_\Omega|u(x)|^p\;dx+\frac{1}{r}\int_\Omega|u(x)|^r\;dx.$$
The same arguments as those used in the case of functional 
$J_\lambda$ show that $I_\lambda$ is well-defined on $E$
and $I_\lambda\in C^1(E,\RR)$ with the the Fr\'echet 
derivative given by
\begin{eqnarray*}
\langle I_\lambda^{'}(u),v\rangle&=&\int_\Omega\log(1+|\nabla 
u(x)|^q)|\nabla u(x)|^{p-2}\nabla u(x)\nabla v(x)\;dx\\
&-&\lambda\int_\Omega|u(x)|^{p-2}u(x)v(x)\;dx+\int_\Omega
|u(x)|^{r-2}u(x)v(x)\;dx\,,
\end{eqnarray*}    
for all $u$, $v\in E$. 
This time our idea is to show that $I_\lambda$ possesses a nontrivial
global minimum point in $E$. We start with the following auxiliary result.

\begin{lemma}\label{PL1}
The functional $I_\lambda$ is coercive on $E$.
\end{lemma}
\proof
In order to prove Lemma \ref{PL1} we first show that for any $b$,
$d>0$ and $0<k<l$ the following inequality holds 
\begin{equation}\label{stea1}
b\cdot t^k-d\cdot t^l\leq b\cdot\left(\frac{b}{d}\right)^{k/(l-k)},
\;\;\;\forall\;t\geq 0.
\end{equation}
Indeed, since the function
$$[0,\infty)\ni t\rightarrow t^\theta$$
is increasing for any $\theta>0$ it follows that
$$b-d\cdot t^{l-k}<0,\;\;\;\forall\;t>\left(\frac{b}{d}\right)^
{1/(l-k)},$$
and
$$t^k\cdot(b-d\cdot t^{l-k})\leq b\cdot t^k< b\cdot\left(\frac{b}{d}
\right)^{k/(l-k)},\;\;\;\forall\;t\in\left[0,\left(\frac{b}{d}
\right)^{1/(l-k)}\right].$$
The above two inequalities show that \eq{stea1} holds true.

Using \eq{stea1} we deduce that for any $x\in\Omega$ and $u\in E$
we have
$$\frac{\lambda}{p}\cdot|u(x)|^p-\frac{1}{r}\cdot|u(x)|^r
\leq\frac{\lambda}{p}\cdot\left[\frac{\lambda\cdot r}{p}\right]
^{(p/(r-p))}=D_1,$$
where $D_1$ is a positive constant independent of $u$ and $x$. 
Integrating the above inequality over $\Omega$ we find
\begin{equation}\label{stea2}
\frac{\lambda}{p}\int_\Omega|u(x)|^p\;dx-\frac{1}{r}
\int_\Omega|u(x)|^r\;dx\leq D_2,\;\;\;\forall\;u\in E
\end{equation}
where $D_2$ is a positive constant independent of $u$.

Using inequalities \eq{l13} and \eq{stea2} we obtain that for
any $u\in E$ with $\|u\|>1$ we have
$$I_\lambda(u)\geq\|u\|^p-D_2.$$
Thus $I_\lambda$ is coercive and the proof of Lemma \ref{PL1} is
complete.  \endproof
{\sc Proof of Theorem \ref{t2}.} 
First, we prove that $I_\lambda$ is weakly lower semicontinuous 
on $E$. Indeed, using the definitions of $J_{\lambda,1}$ and
$J_{\lambda,2}$ introduced in the above section we get
$$I_\lambda(u)=J_{\lambda,1}(u)-J_{\lambda,2}(u),\;\;\;\forall\;
u\in E.$$
Since $\Phi$ is convex it is clear that $J_{\lambda,1}$ is convex
and thus weakly lower semicontinuous on $E$. By Remark \ref{r1} the
functional $J_{\lambda,2}$ is also weakly lower semicontinuous on 
$E$. Thus, we obtain that $I_\lambda$ is weakly lower 
semicontinuous on $E$.

By Lemma \ref{PL1}  we deduce that $I_\lambda$ is coercive 
on $E$. Then Theorem 1.2 in \cite{S} implies that 
there exists $u_\lambda\in E$ a global minimizer of $I_\lambda$ and 
thus a weak solution of problem \eq{2}. 

We show that $u_\lambda$ is not trivial for $\lambda$ large enough.
Indeed, letting $t_0>1$ be a fixed real and $\Omega_1$ be an open 
subset of $\Omega$ with $|\Omega_1|>0$ we deduce that there exists
$u_1\in C_0^\infty(\Omega)\subset E$ such that $u_1(x)=t_0$ for
any $x\in\overline\Omega_1$ and $0\leq u_1(x)\leq t_0$ in $\Omega
\setminus\Omega_1$. We have
\begin{eqnarray*}
I_\lambda(u_1)&=&\int_\Omega\Phi(|\nabla u_1(x)|)\;dx-
\frac{\lambda}{p}\int_\Omega|u_1(x)|^{p}\;dx+\frac{1}{r}
\int_\Omega|u_1(x)|^{r}\;dx\\
&\leq&L-\frac{\lambda}{p}\int_{\Omega_1}|u_1(x)|^{p}\;dx\\
&\leq&L-\frac{\lambda}{p}\cdot t_0^{p}\cdot|\Omega_1|
\end{eqnarray*}
where $L$ is a positive constant.
Thus, there exists $\lambda_\star>0$ such that $I_\lambda(u_1)<0$
for any $\lambda\in[\lambda_\star,\infty)$. It follows that
$I_\lambda(u_\lambda)<0$ for any $\lambda\geq\lambda_\star$ and
thus $u_\lambda$ is a nontrivial weak solution of problem \eq{2}
for $\lambda$ large enough. The proof of Theorem \ref{t2} is 
complete. \endproof

\end{document}